\documentclass[12pt,english]{article}
\usepackage{theorem}
\usepackage{a4wide}
\usepackage{graphics}
\usepackage{epsfig}
\usepackage{amssymb}
\usepackage{latexsym}
\usepackage{amsmath}

\newtheorem{thh}{Theorem}[section]
\newtheorem{df}[thh]{Definition}

\newtheorem{lem}[thh]{Lemma}

\newtheorem{prop}[thh]{Proposition}

\title{Relative cohomology of polynomial mappings}
\author{Philippe Bonnet}

\newcommand{\dem}{{\em Proof: }}
\newcommand{\qed}{\begin{flushright} $\blacksquare$\end{flushright}}

\newcommand{\CC}{\mathbb C}
\newcommand{\FB}{\overline{F}}
\newcommand{\fob}{\overline{f_1}}
\newcommand{\fqb}{\overline{f_q}}
\newcommand{\fib}{\overline{f_i}}
\newcommand{\JT}{J_{\theta}}

\newcommand{\OFB}{\omega_{\overline{F}}}
\newcommand{\CX}{ \mathbb C [x_1,...,x_n]}
\newcommand{\CF}{ \mathbb C [F]}

\newcommand{\CN}{\mathbb C ^n}
\newcommand{\CQ}{\mathbb C ^q}

\begin{document}
\maketitle

\begin{center} { \small
Section de Math\'ematiques, Universit\'e de Gen\`eve, \\ 2-4, rue du li\`evre,
1211 Gen\`eve 24, Switzerland. \\
e-mail: philippe.bonnet@math.unige.ch}
\end{center}

\begin{abstract} Let $F$ be a polynomial mapping from $\CN$ to
$\CQ$ with $n>q$. We study the De Rham cohomology of its fibres
and its relative cohomology groups, by introducing a special fibre
$F^{-1}(\infty)$ "at infinity" and its cohomology. Let us fix a weighted homogeneous
degree on $\CX$ with strictly positive weights. The fibre at
infinity is the zero set of the leading terms of the coordinate
functions of $F$. We introduce the cohomology groups
$H^k(F^{-1}(\infty))$ of $F$ at infinity. These groups enable us to
compute all the other cohomology groups of $F$. For instance, if the
fibre at infinity has an isolated singularity at the origin, we
prove that every weighted homogeneous basis of
$H^{n-q}(F^{-1}(\infty))$ is a basis of all the groups
$H^{n-q}(F^{-1}(y))$ and also a basis of the $(n-q)^{th}$
relative cohomology group of $F$. Moreover the dimension of
$H^{n-q}(F^{-1}(\infty))$ is given by a global Milnor number
of $F$, which only
depends on the leading terms of the coordinate functions of $F$.

\end{abstract}

\section{Introduction}

Let $F$ be a polynomial map from $\CN$ to $\CQ$. In this paper, we
are going to study the De Rham cohomology groups of its fibres and
its relative cohomology groups. These groups have been extensively
studied for holomorphic maps germs, by means of their De Rham
relative complex (see \cite{Loo}, pp. 91 and 137). If $G:(\CN,0)\rightarrow
(\CQ,0)$ denotes a holomorphic map-germ with coordinate functions
$g_1,..,g_q$, and $\Omega^k$ stands for the space of germs of
analytic $k$-forms at the origin of $\CN$, the De Rham relative
complex of $G$ is the following complex: $$ 0 \rightarrow
\CC\{x_1,..,x_n\} \rightarrow \Omega ^1 _G \rightarrow
..\rightarrow \Omega ^k _G \rightarrow \Omega ^{k+1} _G ..$$
defined with the spaces of relative forms: $$ \Omega ^k _G= \Omega
^k/\sum dg_i \wedge \Omega^{k-1} \quad $$ and provided with the
arrows $d_G:\Omega ^k _G\longrightarrow \Omega ^{k+1} _G$ induced
by the exterior derivation. Its cohomology groups $H^k(\Omega ^{*}
_G)$ are called the relative cohomology groups of $G$. They have
been introduced by Hamm and Brieskorn in order to analyse the
topology of isolated singularities defined by holomorphic germs
$g:(\CN,0)\rightarrow (\CC,0)$ (\cite{Ham},\cite{Br}). Their
results have been extended by L\^e D\~ung Tr\'ang, Greuel and Malgrange
to the case of isolated singularities of complete intersections
(\cite{Le},\cite{Gr},\cite{Ma}), and for polynomials $g:\CN
\rightarrow \CC$ satisfying some good conditions at infinity
(\cite{Ga},\cite{B-D}).

Let us go back to the polynomial case, and fix a polynomial map
$F=(f_1,..,f_q)$ from $\CN$ to $\CQ$ where $n>q$. We are going to
construct a special fibre of the polynomial map $F$, that will
play the same role as the singular fibre for holomorphic
map-germs, and define its cohomology. Then we will compute its
cohomology groups and we will see how this cohomology leads us
back to the cohomology of the fibres of $F$.

Let $\Omega^k(\CN)$ be the space of polynomial differential
$k$-forms on $\CN$. By convention we set $\Omega^{-1}(\CN)=0$. For
any ideal $I$ of $\CX$, a polynomial $k$-form $\omega$ is
congruent to zero modulo $I$ ($\omega\equiv 0 \; [I]$) if $\omega$
belongs to $I\Omega^k(\CN)$. In what follows $V(I)$ stands for the
algebraic set of zeros of $I$, i.e : $$V(I)=\{x \in \CN, \;
\forall P \in I, \; P(x)=0\}$$ Recall that the depth of $I$ is the
codimension of $V(I)$, and that $I$ is radical if it is equal to
its root, i.e. to the set of polynomials that vanish on $V(I)$.
Let $F=(f_1,..,f_q)$ be a polynomial map from $\CN$ to $\CQ$ with
$n>q$. We always assume $F$ to be dominant, which means in this
case that its coordinate functions are algebraically independent.
Let us set by convention: $$ \CC[F]=
\CC[f_1,..,f_q]=F^*(\CC[t_1,..,t_q]) $$ Let $deg$ be a weighted
homogeneous degree assigning weights $p_1,..,p_n>0$ to the
variables $x_1,..,x_n$ in $\CX$. Such a degree is called {\em
positive weighted homogeneous} (in short: $p.w.h$). The leading
term of a polynomial $P$ is denoted $\overline{P}$. This degree
can be extended to $\Omega^k (\CN)$ by assigning the degree $p_i$
to $dx_i$: A polynomial $k$-form $\omega=P dx_{i_1} \wedge
..\wedge dx_{i_k}$ is weighted homogeneous of degree $r$ if $P$ is
weighted homogeneous of degree $(r -p_{i_1} -..- p_{i_k})$. By
analogy, we denote by $\overline{\omega}$ the leading term of the
polynomial $k$-form $\omega$. For any point $y=(y_1,..,y_q)$ in
$\CQ$, a polynomial $k$-form $\omega$ is defined to be:

\begin{itemize}
\item{closed on $F^{-1}(y)$ if $d\omega \wedge df_1 \wedge ..\wedge df_q \equiv 0\;[f_1 -y_1,..,f_q-y_q]$,}
\item{exact on $F^{-1}(y)$ if $\omega$ belongs to $d\Omega ^{k-1}(\CN)+(f_1-y_1,..,f_q-y_q)\Omega^k(\CN)$.}
\end{itemize}
Note that $\omega$ is exact on $F^{-1}(y)$ if and only if it belongs to $d\Omega ^{k-1}(\CN)+
\sum \Omega ^{k-1}(\CN)\wedge df_i + (f_1-y_1,..,f_q-y_q)\Omega^k(\CN)$. Indeed, every $k$-form
of type $\eta \wedge df_i$ is exact on $F^{-1}(y)$, as is shown by the formula:
$$
df_i \wedge \eta = d\left \{(f_i -y_i)\eta\right \} - (f_i-y_i)\eta
$$
The $k^{th}$ De Rham cohomology group $H^k(F^{-1}(y))$ is the
quotient of closed $k$-forms on $F^{-1}(y)$ by exact $k$-forms on
$F^{-1}(y)$. These groups are well defined for any point $y$.
Moreover by Grothendieck's theorem (see \cite{Di}, p. 182 or \cite{Gro}), they coincide
with the standard De Rham cohomology groups when $y$ is not a
critical value of $F$. Note that every polynomial $(n-q)$-form is
closed on $F^{-1}(y)$ by definition. From a geometric viewpoint,
we can interpret it by saying that the generic fibres of $F$ have
dimension $(n-q)$. Hence every $(n-q)$-form is closed by restriction
to these fibres. For $k>0$, a polynomial $k$-form $\omega$
is defined to be:

\begin{itemize}
\item{relatively closed if $d\omega \wedge df_1 \wedge ..\wedge df_q =0$,}

\item{relatively exact if $\omega$ belongs to $d\Omega ^{k-1}(\CN) +
 \sum \Omega ^{k-1}(\CN)\wedge df_i$.}
\end{itemize}
The quotient $H^k(F)$ of relatively closed $k$-forms by relatively
exact $k$-forms is the $k^{th}$ relative cohomology group of $F$.
In what follows, we fix a positive weighted homogeneous degree on $\CX$, and
we will not refer to it unless necessary. We denote by $\overline{F}$ the map
$\overline{F}=(\overline{f_1},..,\overline{f_q})$. 

\begin{df}
A polynomial map $F$ is a complete intersection at infinity if the ideal
$I=(\fob,..,\fqb)$ is radical of depth $q$. Its fibre at infinity
is the set $F^{-1}(\infty)=\overline{F}^{-1}(0)$.
\end{df}
Assume the degree $deg$ is the canonical degree on $\CX$. Consider $\CN$ as embedded
in $\mathbb{P}^n(\CC)$ via the map $(x_1,..,x_n)\mapsto [1;x_1;..;x_n]$.
Then $F^{-1}(\infty)$ is the cone in $\CC^{n+1}$ corresponding to the trace
of the fibre $F^{-1}(y)$ on the
hyperplane at infinity. This justifies the terminology of "fibre at infinity".

Let $J$ be the ideal generated by the $q$-minors of the matrix $d\overline{F}$.
The singular set of $F^{-1}(\infty)$ is the set
$Sing(F^{-1}(\infty))=V(I+J)$. Note that $F$ is a complete
intersection at infinity if and only if the ideal $I+J$ has depth
$>q$. A polynomial $k$-form $\omega$ is said to be:
\begin{itemize}
\item{closed at infinity if $d\omega \wedge d\fob \wedge ..\wedge d\fqb
\equiv 0\;[\fob,..,\fqb]$,}
\item{exact at infinity if $\omega$ belongs to $d\Omega ^{k-1}(\CN) +
(\fob,..,\fqb)\Omega^k(\CN)$.}
\end{itemize}
The {\em $k^{th}$ cohomology group at infinity} is the quotient
$H^k(F^{-1}(\infty))$ of closed $k$-forms at infinity by exact
$k$-forms at infinity. Note that by construction $H^k(F^{-1}(\infty))=
H^k(\overline{F}^{-1}(0))$. 

\begin{thh} \label{coho1}
Let $F$ be a complete intersection at infinity. Then
$H^0(F^{-1}(y))=\CC$ for any $y$ in $\CQ$. If $k>0$ and $dim
\;Sing(F^{-1}(\infty))
< n-q-k$, then $H^k(F^{-1}(y))=0$ for any $y$ in $\CQ$.
If $dim \; Sing(F^{-1}(\infty))=0$, then $H^{n-q}(F^{-1}(\infty))$
has dimension $$\mu = dim \; \CX /(J+I)$$ Moreover any weighted
homogeneous basis $\omega_1,..,\omega_{\mu}$ of
$H^{n-q}(F^{-1}(\infty))$ forms a basis of $H^{n-q}(F^{-1}(y))$
for any $y$ in $\CQ$.
\end{thh}
This result is the analogue of what happens in the local case, for singularities
of complete intersections (see \cite{Gr}, pp. 259-260 and 264). If $dim \; Sing(F^{-1}(\infty))=0$,
we can interpret the union of
all the $H^{n-q}(F^{-1}(y))$ as a vector bundle over $\CQ$, whose
space of global sections would be the group $H^{n-q}(F)$. It is
therefore natural to think that these sections are generated by
the base $\omega_1,..,\omega_{\mu}$, and this is exactly what
happens as shown in the following result.

\begin{thh} \label{coho2}
Let $F$ be a complete intersection at infinity. Then $H^0(F)=\CF$.
If $k>0$ and $dim \;Sing(F^{-1}(\infty))< n-q-k$, then $H^k(F)=0$.
If $dim \; Sing(F^{-1}(\infty))=0$, then $H^{n-q}(F)$ is a free
and finitely generated $\CF$-module of rank $\mu$. More precisely
if $\omega_1,..,\omega_{\mu}$ is a weighted homogeneous basis of
$H^{n-q}(F^{-1}(\infty))$, then every polynomial $(n-q)$-form
$\omega$ can be written as: $$ \omega = \sum _i a_i(F)\omega_i +
d\Omega + \sum_j \eta_j \wedge df_j$$ where $a_i$ are uniquely
determined polynomials, and the degrees of the terms of this sum
satisfy the following inequalities: $$ deg(a_i(F)) \leq
deg(\omega) -deg(\omega_i), \quad deg(\Omega)\leq deg(\omega),
\quad deg(\eta_j) \leq deg(\omega) - deg(f_j) $$
\end{thh}
Since $F^{-1}(\infty)$ is the zero set of some weighted homogeneous polynomials, it
is a quasi-cone in $\CN$. Therefore $F^{-1}(\infty)$ has an isolated singularity
at 0 if and only if $dim \; Sing(F^{-1}(\infty))=0$. We then say that {\em $F$ defines
an isolated singularity of complete intersection at infinity}. In this case,
it may be tedious to compute a proper basis of the
cohomology at infinity. The following proposition enables us to
find at least a generating system for it. The contraction along
the Euler vector field $X$ (see \S 2) defines a $\CX$-morphism: $$ i_X:
\Omega ^{n-q}(\CN) \rightarrow \Omega ^{n-q-1}(\CN), \quad \omega
\mapsto i_X(\omega)$$ whose kernel is a noetherian module, hence
finitely generated.

\begin{prop} \label{generators}
Let $F$ be an isolated singularity of complete intersection at infinity. Let
$\{P_k\}$ be a basis of the
algebra $\CX /(I+J)$, and $\{\omega_l\}$ a system of generators of
$ker\; i_X$. Then $\{P_k \omega_l\}_{k,l}$ forms a system of
generators of $H^{n-q}(F^{-1}(\infty))$.
\end{prop}
As a conclusion, we would like to insist on the fact that the
notion of complete intersection at infinity is extrinsic. More
precisely, the fact that an algebraic set $X$ in $\CN$ is given
as a fibre of a complete intersection at infinity strongly
depends on its embedding in $\CN$. For instance there exist 
some plane curves that are not complete intersections
at infinity in $\CC^2$. However they can be embedded in
$\CC^3$ as fibres of complete intersections at infinity.
We illustrate this fact by an example at the end of this paper.
The previous results enable us to compute the cohomology of
some curves by suitably embedding them in affine space, especially if
they do not appear at first sight as
fibres of a complete intersection at infinity. This is quite
surprising since cohomology is completely intrinsic, i.e does not
depend on any particular embedding.

\section{Description of the cohomology at infinity}

In this section, we compute the cohomology groups
$H^{k}(F^{-1}(\infty))$. The main tool will be the
De Rham Lemma on the
division of forms (see \cite{Sai}, p. 166). Let $deg$ be a p.w.h degree.
To that degree corresponds the Euler vector field: $$X=\sum p_i
x_i\frac{\partial }{\partial x_i}$$ and the $\CC ^*$-action
$\varphi$ on $\CN$ defined by
$\varphi_t(x_1,..,x_n)=(t^{p_1}x_1,..,t^{p_n}x_n)$. Let $L_X$ be
the Lie derivative with respect to $X$. As one can easily check, a
polynomial $k$-form $\omega$ is weighted homogeneous of degree $r$
if and only if $\varphi_t ^*(\omega)=t^r \omega$, or equivalently
if $L_X(\omega)= r\omega$.

\begin{prop} \label{infinity}
Let $F$ be a complete intersection at infinity. Then $H^0
(F^{-1}(\infty))=\CC$. If $k>0$ and $dim \;Sing(F^{-1}(\infty))
<n-q-k$, then $H^k(F^{-1}(\infty))=0$. If $dim
\;Sing(F^{-1}(\infty))=0$, then $H^{n-q}(F^{-1}(\infty))$ has
dimension $$\mu = dim \; \CX /(I+J)$$ Moreover if $\{P_k\}$ is a
basis of the algebra $\CX /(I+J)$, and $\{\omega_l\}$ is a system
of generators of $ker\; i_X$, then $\{P_k \omega_l\}_{k,l}$ forms
a system of generators of $H^{n-q}(F^{-1}(\infty))$.
\end{prop}
Note that this result will imply proposition \ref{generators} once
theorem \ref{coho1} has been proved. For simplicity, denote by
$\OFB$ the $q$-form $d\fob \wedge .. \wedge d\fqb$. We first start
by calculating the 0-cohomology, with the following lemma.

\begin{lem} \label{connex}
If $F$ is a complete intersection at infinity, then
$H^0(F^{-1}(\infty))=\CC$.
\end{lem}
\dem Let $R$ be a polynomial such that $dR \wedge \OFB\equiv 0
[I]$. Since $I$ is a radical ideal generated by the $\fib$, the
restriction of $R$ to the fibre at infinity is singular at any
smooth point of this fibre. Hence this restriction is locally
constant. Since the fibre at infinity is defined by weighted
homogeneous polynomials, and $deg$ is a positive degree, this
fibre is a quasi-cone, so it is connected. Therefore $R$ is
constant on the fibre at infinity. There exists a constant
$\lambda$ such that $R-\lambda$ vanishes on $F^{-1}(\infty)$.
Since $I$ is radical, Hilbert's Nullstellensatz implies that
$R\equiv \lambda [I]$, hence giving the result. \qed We pass on to
cohomology of order $k>0$. Let $H^k$ be the following
$\CX$-module:
$$H^k= \frac{\{\omega \in \Omega^k(\CN), \omega
\wedge \OFB \equiv 0\;\;[I]\}} { \sum_i \Omega^{k-1}(\CN)\wedge
d\fib +I \Omega^k(\CN)} $$
We are going to construct a map from
$H^k$ to $H^k(F^{-1}(\infty))$, and see when this leads to an
isomorphism. Let us begin by proving that $H^k$ is annihilated by
$I+J$, a result that is more or less established in a weaker form
in \cite{Sai}, p. 166.

\begin{lem}
For any $k\leq n-q$, we have $(I+J)H^k=0$.
\end{lem}
\dem By definition, $IH^k=0$ and there only remains to check that
$JH^k=0$. Let $\Delta$ be a $q$-minor of the matrix $d\overline{F}$.
By assumption there exist some linear forms $l_{q+1},..,l_{n}$
such that: $$ d\fob \wedge ..\wedge d\fqb \wedge dl_{q+1}\wedge
..\wedge dl_{n}= \Delta dx_1\wedge ..\wedge dx_n$$
Denote the 1-forms $d\fob,..,d\fqb, dl_{q+1},..,dl_{n}$ by $\theta_1,..,\theta_n$. For any polynomial $k$-form
$\omega$ representing an element of $H^k$, there exist some
rational functions $a_{i_1,..,i_k}$ such that: $$\omega =
\sum_{i_1,..,i_k} a_{i_1,..,i_k} \theta_{i_1} \wedge ..\wedge
\theta_{i_k}$$ For any multi-index $(i_1,..,i_k)$, let
$(j_1,..,j_{n-k})$ be a multi-index containing all the elements of
$\{1,..,n\}-\{i_1,..,i_k\}$. Up to a sign $\epsilon$, we get by wedge
product: $$\omega \wedge \theta_{j_1} \wedge ..\wedge
\theta_{j_{n-k}}= \epsilon \Delta a_{i_1,..,i_k} dx_1 \wedge ..\wedge
dx_n$$ Hence $b_{i_1,..,i_k}= \Delta a_{i_1,..,i_k}$ is a
polynomial for any multi-index. Now if the multi-index
$(i_1,..,i_k)$ contains one of the indices $1,..,q$, for instance
$i$, there exists a polynomial $(k-1)$-form $\eta_{i_1,..,i_k}$
such that: $$b_{i_1,..,i_k} \theta_{i_1} \wedge ..\wedge
\theta_{i_k}= \eta_{i_1,..,i_k} \wedge d\fib$$ If the multi-index
$(i_1,..,i_k)$ does not contain any of the indices $1,..,q$, then
all the indices $1,..,q$ are contained in $(j_1,..,j_{n-k})$.
Since $\omega \wedge \omega_{\FB}\equiv 0 \; [I]$, we get: $$ b_{i_1,..,i_k} dx_1 \wedge ..\wedge
dx_n=\omega \wedge \theta_{j_1} \wedge ..\wedge
\theta_{j_{n-k}}\equiv 0\;[I]$$ So $ b_{i_1,..,i_k}$ belongs to
$I$. With these two cases, it is clear that we can find some
polynomial $(k-1)$-forms $\eta_i$ such that: $$ \Delta \omega
\equiv \sum_i \eta_i \wedge d\fib \;[I]$$ So $\Delta H^k=0$ and
the result follows since $J$ is generated by all the $q$-minors
$\Delta$ of $d\FB$. \qed

\begin{lem}
For any $k\leq n-q$, the inclusion map induces a morphism $i_k
^{*} : H^k\rightarrow H^k(F^{-1}(\infty))$.
\end{lem}
\dem Let us show first the inclusion: $$ \{\omega \in
\Omega^k(\CN), \omega \wedge \omega_{\FB}\equiv 0 \;[I]\}
\subset\{\omega \in \Omega^k(\CN), d\omega \wedge
\omega_{\FB}\equiv 0 \;[I]\} $$ Let $\omega$ be a polynomial
$k$-form such that $\omega \wedge\OFB\equiv 0 \;[I]$. Since $I$ is
a complete intersection, $J+I$ has depth $\geq q+1$. Since $\CX$
is catenary, there exists a polynomial $\alpha$ in $J+I$ such that
$\alpha,\fob,..,\fqb$ is a regular sequence. By the previous
lemma, there exist some polynomial forms $\eta_i$ and $\Omega_i$
such that: $$ \alpha\omega = \sum_{i=1} ^q \eta_i \wedge
d\{\overline{f_i}\}+ \sum_{i=1} ^q \overline{f_i}\Omega_i $$ After
derivation and wedge product with $\omega_{\FB}$, we get: $$
\alpha d\omega \wedge \omega_{\FB} + d\alpha\wedge \omega \wedge
\omega_{\FB}\equiv \alpha d\omega \wedge \omega_{\FB}\equiv 0\;[I]
$$ As $\alpha,\fob,..,\fqb$ is a regular sequence, that implies
$d\omega \wedge\omega_{\FB}\equiv 0\;[I]$. Thus the inclusion map
defines a morphism: $$ i_k : \{\omega \in \Omega^k(\CN), \omega
\wedge \omega_{\FB}\equiv 0 \;[I]\} \longrightarrow \{\omega \in
\Omega^k(\CN), d\omega \wedge \omega_{\FB}\equiv 0 \;[I]\} $$
Assume there exist some polynomial forms $\eta_i$ such that
$\omega \equiv \sum_i \eta_i\wedge d\fib\;[I]$. By an easy
computation, we get $\omega\equiv d\{-\sum_i \fib\eta_i\} \;[I]$.
By passage to the quotients, $i_k$ induces a morphism $i_k ^{*} :
H^k\rightarrow H^k(F^{-1}(\infty))$. \qed

\begin{lem} \label{cohinf1}
For any $k\leq n-q$, the morphism $i_k ^*$ is surjective. 
\end{lem}
\dem It is enough to prove surjectivity for weighted homogeneous
$k$-forms. Let $\omega$ be a weighted homogeneous $k$-form of
degree $p$. By an easy computation, we get: $$
i_X(\omega_{\FB})=\sum_i (-1)^{i-1} \overline{f_i} d(\fob)\wedge
..\wedge d(\overline{f_{i-1}})\wedge d(\overline{f_{i+1}})\wedge
..\wedge d(\fqb)\equiv 0\;[I] $$ Starting from the relation
$d\omega \wedge \omega_{\FB}\equiv 0\;[I]$, we obtain after
contraction along the Euler vector field $X$: $$ i_X\left(d\omega
\wedge \omega_{\FB}\right)\equiv i_X(d\omega )\wedge
\omega_{\FB}\equiv 0\;[I] $$ which implies in terms of Lie
derivative: $$ \{L_X(\omega)-d(i_X(\omega))\}\wedge
\omega_{\FB}\equiv 0\;[I] $$ Since $\omega$ is weighted
homogeneous of degree $p$, we have $L_X(\omega)=p\omega$. The
previous congruence can be rewritten as: $$
\{p\omega-d(i_X(\omega))\}\wedge \omega_{\FB}\equiv 0\;[I] $$ Set
$\omega_0 = \omega - d(i_X(\omega)/p)$. By construction,
$\omega_0$ represents an element of $H^k$ for which $i^*
_k(\omega_0)=\omega$. 
\qed

\begin{lem} \label{cohinf2}
Let $F$ be a complete intersection at infinity. If $dim
\;Sing(F^{-1}(\infty)) \leq n-q-k$, then $i_k ^{*}$ is injective.
\end{lem}
\dem Let $\omega$ be a polynomial $k$-form such that $\omega\equiv
d\Omega\; [I]$ and $\omega\wedge \omega_{\FB}\equiv 0\;[I]$. This
implies $d\Omega \wedge \omega_{\FB}\equiv 0\;[I]$. Since
$Sing(F^{-1}(\infty))$ has dimension $\leq n-q-k$, the ideal $J+I$
has codimension $\geq q+k$ in $\CX$. As $I$ is generated by a
regular sequence of length $q$, $\CX / I$ is a Cohen-Macaulay ring
of dimension $(n-q)$. So $J+I$ has depth $\geq k$ in $\CX / I$. By
De Rham lemma ([Sai]), the quotient $H^{k-1}$ is zero. By the
previous lemma, $i_{k-1} ^*$ is surjective, which implies: $$
H^{k-1}(F^{-1}(\infty))=0 $$ There exist some polynomial forms
$\Omega_0$ and $\eta_i$ such that $\Omega=d\Omega_0 + \sum_i \fib
\eta_i$. By an easy computation, we get: $$ \omega\equiv d\Omega
\equiv \eta_1 \wedge d(\fob)+..+\eta_q \wedge d(\fqb)\;[I] $$ Thus
$\omega$ has null class in $H^k$, hence proving injectivity. \qed
{\bf{Proof of proposition \ref{infinity}:}} Let $F$ be a complete
intersection at infinity. Assume that $Sing(F^{-1}(\infty))$ has
dimension $<n-q-k$. For the same reason as in lemma \ref{cohinf2},
the ideal $J+I$ has depth $>k$ in $\CX / I$. By De Rham lemma,
$H^k=0$. Since $i^* _k$ is surjective by lemma \ref{cohinf1},
$H^k(F^{-1}(\infty))=0$.

Assume that $Sing(F^{-1}(\infty))$ has dimension 0. By the
previous lemmas, $i^* _{n-q}$ defines an isomorphism from
$H^{n-q}$ to $H^{n-q}(F^{-1}(\infty))$. By a result due to Greuel
([Gr]), $H^{n-q}$ is finite-dimensional of dimension $$\mu=dim \;
\CX /(J + I)$$ Let $\omega$ be a weighted homogeneous polynomial
$(n-q)$-form of degree $p$ representing an element of $H^{n-q}$.
Following the proof in lemma \ref{cohinf1}, we can see that
$\omega$ is represented by the $(n-q)$-form:
$$
\omega_0 = \omega-d(i_X(\omega/p))
$$
Since $\omega$ is weighted homogeneous of degree $p$, we find:
$$
\omega_0 = \frac{p\omega-d i_X(\omega)}{p}= \frac{L_X(\omega)-d i_X(\omega)}{p}=
\frac{i_X(d\omega)}{p}
$$
Therefore the form $\omega_0$ belongs to $H^{n-q}$ and since $i_X \circ i_X =0$, we
have $i_X(\omega_0)=0$. So every element of $H^{n-q}$ can be represented by
an element of $ker\; i_X$. Conversely, if $\omega_0$ is a weighted homogeneous
$(n-q)$-form of degree $p$ such that $i_X(\omega _0)=0$, then we find:
$$
L_X(\omega_0)=p\omega_0=i_X(d\omega_0)
$$ 
Since $d\omega_0$ is a $(n-q+1)$-form and $\OFB$ is a $q$-form on $\CN$, this
implies:
$$
\omega_0 \wedge \OFB \equiv i_X \left ( \frac{d\omega_0}{p} \right ) \wedge \OFB \equiv
i_X \left (\frac{d\omega_0}{p}\wedge \OFB \right )\equiv 0 \; [I]
$$
So $\omega_0$ belongs to $H^{n-q}$, and $H^{n-q}$ is generated as a $\CX$-module
by the kernel of the contraction morphism.
Let
$\{\omega_l\}$ be a system of generators of $ker \; i_X$. Since
$(I+J)H^{n-q}=0$, $H^{n-q}$ is provided with a structure of $\CX
/(I+J)$-module, with $\{\omega_l\}$ as a generating set. If
$\{P_k\}$ is a basis of the vector space $\CX /(I+J)$, then the
family $\{P_k\omega_l\}_{k,l}$ spans $H^{n-q}$ as a vector space,
hence proving the last assertion of proposition \ref{infinity}.
\qed

\section{The reduction lemma}

In this section, we show how to control the degrees of the
polynomials and polynomial forms occuring in the definition of
closedness and exactness on a fibre of $F$. Beyond the "control"
aspect, this lemma is essential since it will play the role of the
coherence theorems in the local case (see \cite{Loo}, p. 144). In
particular it will enable us to transfer all the cohomological
information from the fibre at infinity to the other fibres of
$F$.\\
\
\\ {\bf{Reduction Lemma}} {\it Let $F$ be a complete intersection
at infinity, $y$ a point in $\CQ$ and $\omega$ a polynomial
$k$-form. If $\omega$ is closed on $F^{-1}(y)$, then
$\overline{\omega}$ is closed at infinity. If $\omega$ is exact on
$F^{-1}(y)$ and $dim \; Sing(F^{-1}(\infty))\leq n-q-k$, then
$\overline{\omega}$ is exact at infinity.} \\ \ \\
The proof proceeds as follows. First we write what it means for a form $\omega$ to be
closed (resp. exact) on a fibre. These definitions involve some polynomials and polynomial
forms. Then we show how to reduce the degrees of these polynomials and polynomial forms whenever
possible. After that, we consider the leading part of $\omega$, and we see that it is closed
(resp. exact) on the fibre at infinity. 

\begin{lem} \label{red1}
Let $f_1,..,f_q$ be a collection of polynomials such that
$\fob,..,\fqb$ is a regular sequence. For any $P$ in $(f_1,..,f_q)$,
there exist some polynomials $a_1,..,a_q$ such that $P=\sum_k a_kf_k$
and $deg(a_i)\leq deg(P)- deg(f_i)$ for any $i$.
\end{lem}
\dem by induction on $q\geq 1$. For $q=1$, this is obvious. Assume
it is true to the order $q$. Let $(f_1,..,f_{q+1})$ be a
collection of polynomials such that $\fob,..,\overline{f_{q+1}}$
is a regular sequence. Let $P$ be an element of
$(f_1,..,f_{q+1})$. Among all the polynomials $a_{q+1}$ such that
$P - a_{q+1}f_{q+1}$ belongs to $(f_1,..,f_q)$, we fix one for
which the degree $r=deg(a_{q+1}f_{q+1})$ is minimal. Let us show
by absurd that $r \leq deg(P)$. Assume that $r >deg(P)$. By the
induction hypothesis, there exist some polynomials $a_i$ such
that: $$ P - a_{q+1}f_{q+1}= \sum_{k=1} ^q a_kf_k \quad \mbox{et}
\quad deg(a_1f_1),..,deg(a_qf_q)\leq r $$ By considering only
terms of degree $r$ in this equality, we get: $$ b_1\fob+..+
b_q\fqb+ b_{q+1}\overline{f_{q+1}}=0 $$ where $b_i$ is either
zero, or the leading term of $a_i$. Since
$\fob,..,\overline{f_{q+1}}$ is a regular sequence,
$\overline{f_{q+1}}$ is not a zero-divisor modulo $\fob,..,\fqb$.
So $b_{q+1}=\overline{a_{q+1}}$ belongs to $(\fob,...,\fqb)$. Let
$\alpha_1,..,\alpha_q$ be some polynomials such that $b_{q+1}=
\sum _i \alpha_i\overline{f_i}$. Since every $\overline{f_i}$ is
weighted homogeneous, we may assume that every $\alpha_i$ is
weighted homogeneous of degree $deg(b_{q+1}) - deg(f_i)$. We set:
$$ c_i = a_i + \alpha_i f_{q+1} \quad \mbox{if} \quad i\not=q+1,
\quad \mbox{and} \quad c_{q+1}= a_{q+1} - \sum_{i=1} ^q
\alpha_i f_i $$ By construction, we deduce $P-
c_{q+1}f_{q+1}= c_1f_1+..+c_q f_q$ and $deg(c_{q+1}f_{q+1})< r$,
hence contradicting the minimality of $r$. Therefore $r\leq
deg(P)$, the polynomial $Q=P - a_{q+1}f_{q+1}$ belongs to
$(f_1,..,f_q)$ and has degree $\leq deg(P)$. Since $\fob,..,\fqb$
is a regular sequence, there exist by induction some polynomials
$a_i$ such that: $$ Q= P - a_{q+1}f_{q+1}=a_1f_1+..+a_q f_q $$ and
whose degrees satisfy the inequalities $deg(a_i)\leq
deg(P)-deg(f_i)$ for any $i$. \qed

\begin{lem} \label{red12}
Let $I'$ be an ideal of $\CX$, $\theta=\{\theta_1,..,\theta_r\}$ a collection of polynomial 1-forms
and $\JT$ the ideal
generated by the $r$-minors of the matrix $(\theta_1,..,\theta_r)$. Let $\eta_1,..,\eta_r$ be some polynomial
$k$-forms such that $\eta_1 \wedge \theta_1+..+\eta_r \wedge \theta_r
\equiv 0\;[I']$. If $\JT + I'$ has depth $>k$ in $\CX /I'$, there exist a
collection $\{\zeta_{i,j}\}$ of polynomial $(k-1)$-forms
such that $\zeta_{i,j}=\zeta_{j,i}$ for all $(i,j)$ and
$\eta_i \equiv \sum_j \zeta_{i,j} \wedge \theta_j \;[I']$ for all $i$.
\end{lem}
\dem Let us show this assertion by induction on $r\geq 1$. For $r=1$, this
is obvious. Indeed, let $\eta_1$ be a polynomial $k$-form such that
$\eta_1\wedge \theta_1\equiv 0 \; [I']$. By assumption, $\JT+I'$ has
depth $>k$ in $\CX /I'$. By De Rham lemma, there exists a polynomial
$(k-1)$-form $\eta_{(1,1)}$ such that $\eta_1\equiv \eta_{(1,1)}
\wedge \theta_1\; [I']$. Assume this assertion holds to the order
$r-1$, and let $\eta=(\eta_1,..,\eta_r)$ be some polynomial $k$-forms
such that:
$$
\eta_1\wedge \theta_1+..+\eta_r\wedge \theta_r\equiv 0\; [I']
$$
By wedge product with $\theta_2,..,\theta_r$, we get:
$$
\eta_1\wedge \theta_1 \wedge ..\wedge \theta_r \equiv 0\; [I']
$$
Since $\JT+I'$ has depth $>k$ in $\CX /I'$, we can apply De Rham lemma.
There exist some polynomial $(k-1)$-forms $\eta_{(1,i)}$ such that:
$$
\eta_1\equiv \eta_{(1,1)}\wedge \theta_1+..+\eta_{(1,r)}\wedge \theta_r\; [I']
$$
For any $i\geq 2$, we set $\overline{\eta_i}=\eta_i - \eta_{(1,i)}
\wedge \theta_1$. By construction, the $k$-forms $\overline{\eta_2}
,..,\overline{\eta_r}$ satisfy the relation:
$$
\sum_{i\geq 2} \overline{\eta_i}\wedge \theta_i\equiv \sum_{i\geq 2}\eta_i\wedge \theta_i +
\sum_{i\geq 2}\eta_{(1,i)}\wedge \theta_i \wedge \theta_1 \equiv \sum_{i\geq 2}\eta_i\wedge \theta_i +
\eta_1 \wedge \theta_1 \equiv 0\; [I']
$$
Let $\overline{\theta}$ be the collection $\{\theta_2,..,\theta_r\}$.
Since the $r$-minors of $(\theta_1,..,\theta_r)$ can be expressed with
the $(r-1)$-minors
of $(\theta_2,..,\theta_r)$, we have the inclusion $\JT
\subset J_{\overline{\theta}}$. So $(J_{\overline{\theta}} + I')$
has depth $>k$ in $\CX /I'$. By the induction's hypothesis, there exist
a collection $(\zeta_{(i,j)})_{2\leq i,j \leq r}$ of polynomial
$(k-1)$-forms, such that $\zeta_{(i,j)}=\zeta_{(j,i)}$ if $2\leq i,j \leq r$
and for which:
$$
\overline{\eta_i} \equiv \zeta_{(i,2)}\wedge \theta_2+..+\zeta_{(i,r)}\wedge
\theta_r \; [I']
$$
This implies for any $i\geq 2$:
$$
\eta_i\equiv \eta_{(1,i)}\wedge \theta_1 +\zeta_{(i,2)}\wedge \theta_2+..
+\zeta_{(i,r)}\wedge \theta_r \; [I']
$$
We extend the collection $(\zeta_{(i,j)})_{2\leq i,j \leq r}$ to a new
collection $(\zeta_{(i,j)})_{1\leq i,j \leq r}$ by setting
$\zeta_{(i,1)}=\zeta_{(1,i)}=\eta_{(1,i)}$ for any $i$.
By construction, $\zeta_{(i,j)}= \zeta_{(j,i)}$ for any $1\leq i,j \leq r$.
Moreover we have for any $i$:
$$
\eta_i\equiv \zeta_{(i,1)}\wedge \theta_1 +\zeta_{(i,2)}\wedge \theta_2+
..+\zeta_{(i,r)}\wedge \theta_r \; [I']
$$
which proves the assertion to the order $r$, and ends this induction.
\qed

\begin{lem} \label{red2}
Let $F$ be a complete intersection at infinity, and an integer
$k>0$ such that $dim \; Sing(F^{-1}(\infty))\leq n-q-k$. Let
$\omega$ be a weighted homogeneous $(k-1)$-form of degree $>0$
such that $d\omega\equiv 0\;[I]$. Then there exists a polynomial
$(k-2)$-form $\Omega$ such that $\omega\equiv d\Omega \; [(I)^2]$.
\end{lem}
\dem We first consider the case $k>1$. Starting from the relation
$d\omega\equiv 0\;[I]$, we get $d\omega\wedge \omega_{\FB} \equiv
0\;[I]$. So $\omega$ is closed at infinity. Since
$Sing(F^{-1}(\infty))$ has dimension $\leq n-q-k$,
$H^{k-1}(F^{-1}(\infty))$ is isomorphic to $H^{k-1}$ by lemmas
\ref{cohinf1} and \ref{cohinf2}. By De Rham lemma, these quotients
are zero. So $\omega$ is exact at infinity. There exist some
polynomial forms $\Omega$, $\eta_i$ such that: $$ \omega = d\Omega
+\fob \eta_1+..+\fqb \eta_q $$ This yields after derivation: $$
d\omega\equiv \eta_1 \wedge d\fob+..+\eta_q \wedge d\fqb \equiv 0
\;[I] $$ Since $dim \;Sing(F^{-1}(\infty))\leq n-q-k$, the ideal
$J+I$ has depth $\geq q+k$ in $\CX$. Since $I$ is generated by a
regular sequence, $J+I$ has depth $\geq k$ in $\CX /I$. By lemma
\ref{red12}, there exists a collection $(\zeta_{(i,j)})$ of
polynomial $(k-2)$-forms such that $\zeta_{(i,j)}=\zeta_{(j,i)}$
for any $(i,j)$ and for which: $$ \eta_i\equiv \sum_{j=1} ^q
\zeta_{(i,j)}\wedge d\fib  \;[I] $$ By combining these
congruences, we find: $$ \omega \equiv d\Omega + \sum_{i,j=1} ^q
\zeta_{(i,j)} \wedge \overline{f_i} d(\overline{f_j}) \; [(I)^2]
$$ Since the collection $(\zeta_{(i,j)})$ is symmetric, we can
rewrite it as follows: $$ \omega \equiv d\Omega + \sum_{i<j}
\zeta_{(i,j)} \wedge \{\overline{f_i} d(\overline{f_j})
+\overline{f_j} d(\overline{f_i})\} + \sum_i \zeta_{(i,i)}\wedge
\overline{f_i} d(\overline{f_i}) \; [(I)^2] $$ Let us set: $$
\eta=\sum_{i<j}\overline{f_i}\overline{f_i}\eta_{(i,j)}+ \sum_i
(\overline{f_i})^2 \eta_{(i,i)}/2 $$ By an integration by parts,
we deduce $\omega \equiv d\Omega +d\eta \; [(I)^2]$.

If now $k=1$,
consider a polynomial 0-form $\omega$ such that $d\omega\equiv
0[I]$. By lemma \ref{connex}, there exist a constant $\lambda$ and
some polynomials $a_k$ such that: $$\omega = \lambda + \sum_i a_i
\fib$$ Since $\omega$ is weighted homogeneous of degree $>0$, we
have $\lambda=0$. By derivation we get: $$ d\omega \equiv \sum_i
a_i d\fib \equiv 0 \; [I]$$ By wedge product we find $a_i \OFB
\equiv 0\;[I]$. Since $I+J$ has depth $>q$, this yields $a_i\equiv
0\; [I]$ and $\omega \equiv 0\;[I^2]$. \qed

\begin{lem}
Let $F$ be a complete intersection at infinity. Let
$y=(y_1,..,y_q)$ be a point in $\CQ$, and let $\omega$ be an exact
$k$-form on $F^{-1}(y)$. If $dim \; Sing(F^{-1}(\infty))\leq
n-q-k$, there exist a polynomial $(k-1)$-form $\Omega$ and some
polynomial $k$-forms $\eta_i$ such that: $$ \omega=d\Omega +\sum_i
(f_i-y_i)\eta_i \quad \mbox{and} \quad
deg(\Omega),deg(f_1\eta_1),.., deg(f_q\eta_q)\leq deg(\omega) $$
\end{lem}
\dem The case $k=0$ has already been treated in lemma \ref{red1},
so we pass on to the case $k>0$. Let $\omega$ be an exact $k$-form
on $F^{-1}(y)$. Among all the $(k-1)$-forms $\Omega$ such that
$\omega\equiv d\Omega\;[I]$, let us fix one of minimal degree
$r=deg(\Omega)$. Let us show by absurd that $r\leq deg(\omega)$.
Assume that $r>deg(\omega)$. The form $\omega -d\Omega$ has degree
$r$ and all its coefficients belong to the ideal
$(f_1-y_1,..,f_q-y_q)$. By applying lemma \ref{red1} to these
coefficients, we can see there exist some polynomial $k$-forms
$\eta_i$ such that: $$ \omega=d\Omega +\sum (f_i-y_i)\eta_i \quad
\mbox{and} \quad deg(f_1\eta_1),..,deg(f_q\eta_q)\leq r $$ By
considering only terms of degree $r$ in this equality, we find: $$
d\overline{\Omega}+ \sum \overline{f_i}\zeta_i=0 $$ where
$\zeta_i$ is either zero or the leading term of $\eta_i$. By lemma
\ref{red2}, there exist some polynomial forms $\zeta_{i,j}$ such
that: $$ \overline{\Omega}=\alpha + \sum_{i,j}
\overline{f_i}\overline{f_j}\zeta_{i,j} $$ where $\alpha$ is
either a constant, or an exact form. Since every $\overline{f_i}$
is weighted homogeneous, we may assume that $\alpha$ (resp.
$\zeta_{i,j}$) is weighted homogeneous of degree $r$ (resp.
$r-deg(f_i)-deg(f_j)$). Let us set: $$ \Omega'=\Omega - \alpha -
\sum_{i,j} (f_i-y_i)(f_j-y_j)\zeta_{i,j} $$ By construction, we
have $\omega\equiv d\Omega\equiv d\Omega' [f_1-y_1,..,f_q-y_q]$
and $deg(\Omega')<r$, hence contradicting the minimality of $r$.
So $r\leq deg(\omega)$. By applying lemma \ref{red1} to the
coefficients of $\omega -d\Omega$, we can see there exist some
polynomial $k$-forms $\eta_i$ such that: $$ \omega=d\Omega +\sum
(f_i-y_i)\eta_i \quad \mbox{and} \quad deg(f_1\eta_1),..,
deg(f_q\eta_q)\leq deg(\omega) $$ \qed {\bf{Proof of the reduction
lemma:}} Let $\omega$ be a closed $k$-form on $F^{-1}(y)$, that is
$d\omega \wedge \omega_F\equiv 0 \; [f_1-y_1,..,f_q-y_q]$. By
applying lemma \ref{red1} to the coefficients of $d\omega \wedge
\omega_F$, we can see there exist some polynomial $(k+q)$-forms
$\Omega_i$ such that: $$ d\omega \wedge \omega_F=\sum_k
(f_k-y_k)\Omega_k \quad \mbox{and} \quad
deg(f_1\Omega_1),..,deg(f_q\Omega_q)\leq deg(d\omega \wedge
\omega_F) $$ By considering only terms of degree
$r=deg(\omega)+deg(f_1..f_q)$ in this equality, we find: $$
d(\overline{\omega}) \wedge \omega_{\FB}=\fob\zeta_1+..+\fqb
\zeta_q $$ where $\zeta_i$ is either zero or the leading term of
$\Omega_i$. So $d(\overline{\omega}) \wedge \omega_{\FB}\equiv
0\;[I]$ and $\overline{\omega}$ is closed at infinity. Assume now
that $\omega$ is exact on $F^{-1}(y)$, that $k>0$ and that
$Sing(F^{-1}(\infty))$ has dimension $\leq n-q-k$. By lemma
\ref{red2}, there exist some polynomial $k$-forms $\Omega$ and
$\eta_i$ such that: $$ \omega=d\Omega +\sum_i (f_i-y_i)\eta_i
\quad \mbox{and} \quad deg(\Omega),deg(f_1\eta_1),..,
deg(f_q\eta_q)\leq deg(\omega) $$ By considering only terms of
degree $r=deg(\omega)$ in this equality, we find: $$
\overline{\omega}=d\eta +\sum_i \overline{f_i}\zeta_i $$ where
$\eta$ (resp. $\zeta_i$) is either zero or the leading term of
$\Omega$ (resp. $\eta_i$). Therefore $\overline{\omega}$ is exact
at infinity. The case $k=0$ is treated in exactly the same way.
\qed

\section{Proof of theorem \ref{coho1}}

Let $F$ be a complete intersection at infinity for a positive
weighted homogeneous degree. In this section, we are going to
establish theorem \ref{coho1}. We will split the proof in three
steps.
\\ \ \\ {\bf Assertion 1:} {\it $H^{0}(F^{-1}(y))=\CC$ for any $y$.} \\ \
\\
Let $y=(y_1,..,y_q)$ be a point in $\CQ$. Let us prove by
induction on $r\geq 0$ that for every polynomial $R$ of degree
$\leq r$, closed on $F^{-1}(y)$, there exists a constant $\lambda$
such that: $$R\equiv \lambda \;[f_1 - y_1,..,f_q - y_q]$$ For
$r=0$, this is obvious since every polynomial of degree 0 is a
constant. Assume this is true to the order $(r-1)$. Let $R$ be a
polynomial of degree $\leq r$ closed on $F^{-1}(y)$. By the
reduction lemma, $\overline{R}$ is closed at infinity. By
lemma \ref{connex}, there exist a constant $\lambda$ and some
polynomials $a_i$ such that: $$\overline{R}=\lambda + \sum_i a_i
\fib$$ Since $\overline{R},\fob,.., \fqb$ are weighted
homogeneous, we may assume that every $a_i$ is weighted
homogeneous of degree $r-deg(f_i)$, and that $\lambda=0$. Set
$$R'= R-\sum_i a_if_i$$ By construction $deg(R')<r$ and $R'$ is
closed on $F^{-1}(y)$. Thus there exists a constant $\lambda$ such
that: $$ R\equiv R' \equiv \lambda \;[f_1 - y_1,..,f_q - y_q]$$
hence proving our induction. So the constant function 1 spans the
vector space $H^0(F^{-1}(y))$, and there remains to check that
$(f_1 -y_1,..,f_q-y_q) \not=(1)$. Assume this is not true, and
take some polynomials $a_i$ such that: $$ a_1 ( f_1- y_1)+
...+a_q(f_q-y_q)=1$$ By lemma \ref{red1}, we may assume that every
$a_i$ has degree $\leq -deg(f_i)$, which is obviously impossible.
\\
\
\\ {\bf Assertion 2:} {\it If $dim \; Sing(F^{-1}(\infty))<n-q-k$,
then $H^{k}(F^{-1}(y))=0$ for any $y$.}\\ \ \\ Let us show by
induction on $r$ that every $k$-form $\omega$ of degree $r$,
closed on $F^{-1}(y)$ is exact on $F^{-1}(y)$. This is clear for
$r=min\{p_{i_1}+..+p_{i_k}\}$ because every $k$-form with this
degree has constant coefficients, so it is exact. Assume this
assertion holds to the order $(r-1)$. Let $\omega$ be a polynomial
$k$-form of degree $r$, closed on $F^{-1}(y)$. By the reduction
lemma, $\overline{\omega}$ is closed at infinity. By proposition
\ref{infinity}, $\overline{\omega}$ is exact at infinity
because $Sing(F^{-1}(\infty))$ has dimension $<n-q-k$. There exist
some weighted homogeneous forms $\Omega$ and $\eta_i$ such that:
$$ \overline{\omega}=d\Omega+\fob\eta_1+..+\fqb\eta_q $$ By
construction, $\Omega$ (resp. $\eta_i$) is either zero or has
degree $deg(\omega)$ (resp. $deg(\omega)-deg(f_i)$). The $k$-form
$\omega'=\omega - d\Omega - \sum_i (f_i-y_i)\eta_i$ is closed on
$F^{-1}(y)$ and has degree $\leq (r-1)$. By the induction
hypothesis, $\omega'$ is exact on $F^{-1}(y)$. So $\omega$ is
exact on $F^{-1}(y)$, which proves our induction.\\ \ \\
 {\bf
Assertion 3:} {\it If $dim \; Sing(F^{-1}(\infty))=0$, then every
weighted homogeneous basis $\omega_1,.., \omega_{\mu}$ of
$H^{n-q}(F^{-1}(\infty))$ forms a basis of $H^{n-q}(F^{-1}(y))$
for any $y$.}\\
\
\\
By proposition \ref{infinity}, $H^{n-q}(F^{-1}(\infty))$ is
finite-dimensional of dimension: $$\mu=dim \; \CX / (J+I)$$ Let
${\cal{B}}=\{\omega_1,.., \omega_{\mu}\}$ be a set of weighted
homogeneous $(n-q)$-forms that gives a basis of
$H^{n-q}(F^{-1}(\infty))$. Let us show that ${\cal{B}}$ is a basis
of all the groups $H^{n-q}(F^{-1}(y))$.

Let us prove by absurd that ${\cal{B}}$ is linearly independent in
$H^{n-q}(F^{-1}(y))$. Assume there exist some constants
$a_1,..,a_{\mu}$, not all zero, such that
$\omega=a_1\omega_1+..+a_{\mu}\omega_{\mu}$
is exact on $F^{-1}(y)$. Let $r$ be the maximum of the degrees of the
$\omega_i$ for which $a_i$ is not zero. Since $\omega_1,
..,\omega_{\mu}$ are weighted homogeneous and linearly independent
in $\Omega^{n-q}(\CN)$, $r$ is equal to the degree of $\omega$.
By the reduction lemma, $\overline{\omega}$ is exact at infinity,
which means:
$$
\sum_{deg(\omega_i)=r} a_i \omega_i = 0 \quad \mbox{in} \quad H^{n-q}(F^{-1}(\infty))
$$
Therefore, $a_i=0$ if $deg(\omega_i)=r$, hence contradicting the definition
of $r$. So $\omega_1,..,\omega_{\mu}$ are linearly independent
in $H^{n-q}(F^{-1}(y))$ for any $y$.

Let us show by induction on $r$ that every $(n-q)$-form of degree $r$
is spanned by ${\cal{B}}$ in $H^{n-q}(F^{-1}(y))$. For
$r=min\{p_{i_1}+..+p_{i_{n-q}}\}$, this is clear because
every $(n-q)$-form with this degree has constant coefficients,
so it is exact. Assume this assertion holds to the order
$(r-1)$. Let $\omega$ be a polynomial $(n-q)$-form of degree $r$.
By definition $\overline{\omega}$ is closed at infinity.
So it can be expanded as follows:
$$
\overline{\omega}=\sum_{k=1} ^{\mu} a_k\omega_k+ d\Omega +
\sum_{i=1} ^q \fib \eta_i
$$
where $a_i$ are constant and equal to zero if $deg(\omega_i)\not=r$,
and $\Omega$ (resp. $\eta_i$) is weighted homogeneous de degree $r$
(resp. $r-deg(f_i)$). Let us set:
$$
\omega'=\omega - \sum_{k=1} ^{\mu} a_k\omega_k - d\Omega - \sum_{i=1} ^q
(f_i-y_i)\eta_i
$$
By construction, $\omega'$ has degree $<r$. By the induction hypothesis,
$\omega'$ is spanned by $\omega_1,..,\omega_{\mu}$ in $H^{n-q}(F^{-1}(y))$.
Therefore $\omega$ is also spanned by $\omega_1,..,\omega_{\mu}$ in
$H^{n-q}(F^{-1}(y))$, which proves our induction.
\qed

\section{Proof of theorem \ref{coho2}}

In this section we are going to prove separately the three
assertions of theorem \ref{coho2}, by merely using the same
methods as in theorem \ref{coho1}.

\subsection{Relative $0$-cohomology}

Let $F$ be a complete intersection at infinity, and let us show
that $H^0(F)$ is equal to $\CF$. Consider a polynomial $R$ such
that: $$dR \wedge df_1 \wedge ..\wedge df_q =0$$ Then $R$ is
closed on every fibre of $F$. By theorem \ref{coho1}, we have
$H^0(F^{-1}(y))=\CC$ for any $y$, and $R$ is constant on every
fibre of $F$. Consider the map: $$\alpha: \CQ \longrightarrow \CC,
\quad y\longmapsto \mbox{"unique value of $R$ along
$F^{-1}(y)$"}$$ Its graph corresponds to the image of the mapping
$(f_1,..,f_q,R)$, hence it is a constructible set whose closure is
irreducible. So $\alpha$ defines a rational correspondence in the
sense of Zariski. By Zariski's Main Theorem (see \cite{Mu}, p. 52),
$\alpha$ is a rational function on $\CQ$. Therefore $R$ can be
written as $R=\alpha(F)=A(F)/B(F)$, where $A$ and $B$ are
relatively prime polynomials.

Let us show by absurd that $B$ is a non-zero constant. Assume this
is not. For any point $y$, the fibre $F^{-1}(y)$ is non-empty
because $H^0(F^{-1}(y))=\CC$ by theorem \ref{coho1}. For any point
$y$ in $B ^{-1}(0)$, there exists a point $x$ such that $F(x)=y$,
and so $B(y)R(x)=A(y)=0$. Thus $A$ vanishes on the hypersurface $B
^{-1}(0)$. By Hilbert's Nullstellensatz, $A$ and $B$ cannot be
relatively prime, hence a contradiction.

\subsection{Relative $k$-cohomology}

Let $F$ be a complete intersection at infinity and $k$ an integer $>0$ such that
$Sing(F^{-1}(\infty))$ has dimension $<n-q-k$. We are going to
prove that $H^k(F)=0$.

\begin{lem} \label{first1}
Let $F$ be a complete intersection at infinity such that
$Sing(F^{-1}(\infty))$ has dimension $<n-q-k$. Then the ideal $J$
has depth $>k+1$.
\end{lem}
\dem Let $\FB$ be the map from $\CN$ to $\CQ$ defined in the introduction.
By construction,
$S=V(J)$ is the singular set of $\FB$. This set is globally invariant with respect to the
$\CC ^*$-action $\varphi$, because $J$ is generated by weighted
homogeneous polynomials. Since the weights are all strictly
positive, every irreducible component of $S$ contains the origin
in $\CN$. By the generic smoothness theorem, the closure $Y$ of
$\FB(S)$ has dimension $<q$. By applying the theorem on the
dimension of fibres to $\FB:S\rightarrow Y$, we find: $$ dim \;
Sing(F^{-1}(\infty)) = dim \; S\cap \FB ^{-1}(0) \geq dim \; S -
dim \; Y $$ So $S$ has dimension $<n-k-1$. Since $\CX$ is
catenary, $J$ has depth $>k+1$. \qed

\begin{lem} \label{first2}
Let $F$ be a polynomial map for which the ideal $J$ has depth
$>k+1$. A weighted homogeneous $k$-form $\omega$ satisfies the
equation $d\omega\wedge \omega_{\FB}=0$ if and only if there exist
some weighted homogeneous forms $\Omega$ and $\eta_i$ such that
$\omega = d\Omega + \sum \eta_i \wedge d\fib$
\end{lem}
\dem Let $\omega$ be a weighted homogeneous $k$-form $\omega$, of
degree $r$, such that $d\omega\wedge \omega_{\FB}=0$. By De Rham
lemma, there exists some polynomial $(k-1)$-forms $\zeta_i$ such
that: $$ d\omega= \sum_{i=1} ^q \zeta_i \wedge d\fib $$ So
$\omega$ is a closed form of the relative De Rham complex of
$\FB$. By a result of Malgrange (\cite{Ma}, p. 68), there exist some
germs of analytic forms $\Omega'$ and $\eta_i '$, defined in a
neighborhood of 0 and such that: $$ \omega= d\Omega' + \sum \eta_i
' \wedge d\fib $$ Let $\Omega$ (resp. $\eta_i$) be the weighted
homogeneous part of $\Omega '$ (resp. $\eta_i '$) of degree $r$
(resp. $r-deg(f_i)$). These forms are all polynomials. Since the
forms $\omega$ and $d\fib$ are weighted homogeneous, we get the
equality: $$ \omega= d\Omega + \sum \eta_i  \wedge d\fib $$ \qed
{\bf{Proof of the first part of theorem \ref{coho2}:}} Let $F$ be
a complete intersection at infinity and $k$ an integer $>0$ such that
$Sing(F^{-1}(\infty))$ has dimension $<n-q-k$. Let us prove by
induction on $r$ that every relatively closed $k$-form $\omega$ of
degree $r$ is relatively exact. For $r=1$, this is clear because
every $k$-form of degree 1 has constant coefficients, so it is
exact. Assume this assertion holds to the order $(r-1)$. Let
$\omega$ be a relatively closed $k$-form of degree $r$. By
definition, it satisfies the equation: $$ d\omega \wedge df_1
\wedge ..\wedge df_q=0 $$ By considering only terms of degree
$r+deg(f_1..f_q)$ in this equality, we find
$d(\overline{\omega})\wedge \omega_{\FB}=0$. By lemmas \ref{first1}
and \ref{first2}, there exist some weighted homgeneous forms $\Omega$ and
$\eta_i$ such that: $$ \overline{\omega}= d\Omega + \eta_1 \wedge
d\fob+..+\eta_q \wedge d\fqb $$ Let us set: $$ \omega ' = \omega -
d\Omega - \eta_1 \wedge df_1-..-\eta_q \wedge df_q $$ By
construction $\omega'$ is relatively closed of degree $<r$. By the
induction hypothesis, $\omega'$ is relatively exact. So $\omega$
is relatively exact, which proves our induction. Therefore
$H^k(F)$ is zero. \qed

\subsection{Relative $(n-q)$-cohomology}

Let $F$ be a complete intersection at infinity such that
$Sing(F^{-1}(\infty))$ has dimension zero. We are going to prove
that $H^{n-q}(F)$ is a free and finitely generated module of rank
$\mu$. More precisely fix a weighted homogeneous basis
$\omega_1,..,\omega_{\mu}$ of $H^{n-q}(F^{-1}(\infty))$. We are
going to show that every $(n-q)$-form $\omega$ can be written as:
$$ \omega = \sum_{i=1} ^{\mu} a_i(F)\omega_i + d\Omega +
\sum_{k=1} ^q \eta_k \wedge df_k $$ where the polynomials $a_i$
are uniquely determined, and the degrees of the terms of this sum
satisfy the following inequalities: $$ deg(a_i(F)) \leq
deg(\omega) -deg(\omega_i), \quad deg(\Omega)\leq deg(\omega),
\quad deg(\eta_i) \leq deg(\omega) - deg(f_i) $$ Let us first
prove the existence of such a decomposition, by induction on the
degree $r$ of the polynomial $(n-q)$-form $\omega$. For $r=1$, this assertion
is clear because every $(n-q)$-form of degree 1 is exact. Assume
this holds to the order $(r-1)$. Let $\omega$ be a polynomial
$(n-q)$-form of degree $r$. By theorem \ref{coho1}, there exist
some constants $\lambda_1,.., \lambda_{\mu}$ such that: $$ \omega
'=\omega - \lambda_1\omega_1 - .. - \lambda_{\mu}\omega_{\mu} $$
is exact on $F^{-1}(0)$. Let us show by absurd that $\lambda_i=0$
if $deg(\omega_i)>r$. Assume there exists an index $i_0$ for which
$\lambda_{i_0}\not=0$ and $deg(\omega_{i_0})>p$. Let $p$ be the
maximum of the degrees of the forms $\omega_i$ such that
$\lambda_i\not=0$. Since the $\omega_i$ are weighted homogeneous
and linearly independent in $\Omega ^{n-q}(\CN)$, $p$ is equal to
the degree of $\omega '$. By the reduction lemma, the form: $$
\overline{\omega '}=-\sum_{deg(\omega_i)=p} \lambda_i \omega_i $$
is exact at infinity. So $\lambda_i =0$ if $deg(\omega_i)=p$,
hence contradicting the definition of $p$.

By construction,
$\omega'$ is exact on $F^{-1}(0)$ and has degree $\leq r$. By the
reduction lemma, there exist some polynomial forms
$\Omega,\eta_1,..,\eta_q$ such that: $$ \omega '= \omega -
\sum_{i=1} ^{\mu} \lambda_i\omega_i= d\Omega + \sum_{i=1} ^q
f_i\eta_i $$ where $\Omega$ (resp. $\eta_i$) has degree $\leq r$
(resp. $\leq r-deg(f_i)$). By construction, every $\eta_i$ has
degree $<r$. By the induction hypothesis, there exist some
polynomials $b_{(i,j)}$ and some polynomial forms
$\Omega_i,\zeta_{(i,j)}$ such that: $$ \eta_i= \sum_{j=1} ^{\mu}
b_{(i,j)}(F)\omega_j + d\Omega_i +\sum_{j=1} ^q \zeta_{(i,j)}
\wedge df_j $$ and whose degrees satisfy the following
inequalities: $$ deg(b_{(i,j)}(F)) \leq deg(\eta_i)
-deg(\omega_j), \quad deg(\Omega_i)\leq deg(\eta_i), \quad
deg(\zeta_{(i,j)}) \leq deg(\eta_i) - deg(f_j) $$ From that, we
deduce: $$ \omega= \sum_{i=1} ^{\mu} \left (\lambda_i + \sum_{j=1}
^{\mu} b_{(i,j)}(F)\right )\omega_i + d\left (\sum_{i=1} ^q
f_i\Omega_i \right ) + \sum_{j=1} ^q \left (\sum_{i=1} ^q f_i
\zeta_{(i,j)}- \Omega_j \right ) \wedge df_j $$ It is
straightforward to check the degrees inequalities for all the
terms of this sum. This proves our induction and the existence of
this decomposition. To prove the uniqueness of the polynomials
$a_i$, assume that: $$ \sum_{i=1} ^{\mu} a_i(F)\omega_i + d\Omega
+ \eta_1 \wedge df_1 +.. +\eta_q \wedge df_q =0 $$ Starting from
the equality $\eta_i\wedge df_i=d((f_i-y_i)\eta_i) -
(f_i-y_i)d\eta_i$, we get that $\sum_i a_i(y)\omega_i$ is exact on
$F^{-1}(y)$, for any $y$ in $\CQ$. By theorem \ref{coho1},
$a_i(y)=0$ for any $i$ and any $y$. So all $a_i$ are zero.
Therefore the $\omega_i$ form a basis of $H^{n-q}(F)$. \qed

\section{An example}

We are going to consider an example of polynomial in two
variables, which does not define a complete intersection at
infinity in $\CC^2$. However its generic fibres can be embedded in
$\CC^3$ so as to correspond to the fibres of a complete
intersection at infinity. This enables us to compute the
cohomology of these fibres. Nevertheless we do not know whether, given
an affine curve $C$, it is always possible to embed it into an affine
space as the fibre of a complete intersection at infinity.

Let $f$ be the polynomial $x^4 + x^2y^2$. Obviously, $f$ is not
semi-weighted homogeneous for any degree because $x^2$ divides $f$. For $\lambda\not=0$, the
polynomial $x$ is invertible on the fibre $f^{-1}(\lambda)$. We
introduce a new variable $z$. Then $f^{-1}(\lambda)$ is isomorphic to the
curve in $\CC^3$ given by the equations: $$ xz=1, \quad x^2 + y^2
- \lambda z^2 =0 $$ Let $deg$ be the standard homogeneous degree
on $\CC[x,y,z]$. For $\lambda\not=0$, the map $F(x,y,z)=(xz,x^2 +
y^2 - \lambda z^2)$ defines a complete intersection at infinity.
To see this, it is enough to check that $I + J$ has finite
codimension in $\CC[x,y,z]$. By an easy computation, we find $I +
J=(xy,yz,xz,x^2+ y^2 - \lambda z^2,x^2+ \lambda z^2)$. Therefore
$\{1,x,y,z,x^2\}$ is a basis of the algebra $\CC[x,y,z]/(I + J)$
and $\mu= 5$. By using proposition \ref{generators}, and after
some computations, we obtain the following basis of
$H^1(F^{-1}(\infty))$ for $\lambda\not=0$: $$ \omega_1 = zdx
-xdz,\quad \omega_2= ydz - zdy,\quad \omega_3= xdy -ydx, \quad
\omega_4= x\omega_2, \quad \omega_5= z\omega_1 $$ We easily check that, via
the previous embedding, the projective closure of
$f^{-1}(\lambda)$ in ${\CC} {\mathbb{P}}^ 3$ meets transversally four times the
hyperplane at infinity. Therefore
$f^{-1}(\lambda)$ is a torus that has been punctured four times.

\end{document}